\documentclass[a4paper,
               headexclude,
               headsepline,
               footexclude,
               footsepline,
               plainfootsepline,
               abstracton,
               twoside,
               11pt]{scrartcl}

\usepackage[active]{srcltx} 
\usepackage{scrlayer-scrpage}
\usepackage{amsmath,amsfonts,amsthm,amssymb}
\usepackage{hyperref} 
\usepackage{xcolor}%
\usepackage{charter}
\usepackage[T1]{fontenc}
\usepackage{tikz}
\newcommand{\coloneqq}{\mathrel{\mathop:}=}
\newcommand{\set}[2]{\left\{{#1}\left|\vphantom{#1#2\strut}\right.\, 
                    {#2}\right\}}

\newcommand{\Aa}{\mathbb A}
\newcommand{\DD}{\mathbb D}
\newcommand{\PP}{\mathbb P}
\newcommand{\UU}{\mathbb U}
\newcommand{\XX}{\mathbb X}
\newcommand{\cB}{\mathcal{B}}
\newcommand{\cC}{\mathcal{C}}
\newcommand{\cE}{\mathcal{E}}
\newcommand{\cL}{\mathcal{L}}
\newcommand{\cM}{\mathcal{M}}
\newcommand{\cX}{\mathcal{X}}
\newcommand{\CG}[1]{\Gamma_{#1}}%% confluence graph
\addtokomafont{title}{\rmfamily}
\addtokomafont{section}{\rmfamily}
\addtokomafont{minisec}{\normalfont\rmfamily}%% for "abstract"
\swapnumbers
\theoremstyle{plain}
\newtheorem{theo}{Theorem}[section]
\newtheorem*{main}{Main Result}
\newtheorem{lemm}[theo]{Lemma}
\newtheorem{coro}[theo]{Corollary}
\newtheorem{prop}[theo]{Proposition}
\theoremstyle{definition}

\newtheorem{exam}[theo]{Example}
\newtheorem{exas}[theo]{Examples}
\newtheorem{rema}[theo]{Remark}

\newtheorem{assu}[theo]{Assumption}
\title{Confluence Graphs of Unitals}

\author{Theo Grundh\"ofer, Markus J. Stroppel, Hendrik Van Maldeghem}
\begin{document}
\maketitle

\begin{abstract}
  \noindent%
    We show that the cliques of maximal size in the confluence graph
    of an arbitrary unital of order $q>2$ have size $q^2$, and that
    these cliques are the pencils of all blocks through a given point.
    We also determine all maximal cliques of the confluence graph of
    the Hermitian unitals. %
    As an application, we show that the confluence graph of an
    arbitrary unital unambiguously determines the unital. Along the
    way, we show that each linear space with $q^2$ points such that
    the sizes of both point rows and line pencils are bounded above by
    $q+1$ embeds in a projective plane of order~$q$.
\end{abstract}

\section*{Introduction}

For any incidence geometry $\XX = (X,\cB,I)$, the confluence graph
$\CG\XX$ has as vertices the set $\cB$ of blocks of~$\XX$, and two of
these vertices form an edge if they share a point in~$\XX$. While this
graph is rather boring for projective or affine planes, it sometimes
allows to obtain a fresh look at certain incidence geometries (e.g.,
the confluence graph of the inversive plane of order~$3$ is the
opposite graph of the incidence graph of the generalized quadrangle of
order~$2$, see~\cite{MR3135960}). %
For unitals of order $q>2$, we obtain that the confluence graph
contains complete information about the unital (see
Theorem~\ref{CliquesOfSizeQQ-EKR} below): %

\begin{main}
  Any set of size $q^2$ consisting of mutually intersecting blocks in
  a unital of order $q>2$ coincides with the set of all blocks through
  some point. %
\end{main}

The heart of our argument is an embedding result for certain linear
spaces, see Theorem~\ref{allFromProjective} below. %
As a corollary, we obtain (see Theorem~\ref{reconstructUnital}
below) %
that a unital is completely determined by its confluence graph (whose
vertices are the blocks, adjacent if they share a point).

The classical Erd\H os-Ko-Rado problem asks for the maximal number of
mutually intersecting $k$-subsets of an $n$-set and also asks to
characterize the sets attaining the upper bound. This has been taken
over by geometers, first as a $q$-analogue in vector spaces and
projective spaces, then exploring further analogies for different
geometries by, for instance, asking for the maximal number of mutually
intersecting subspaces of a certain kind in a given geometry. One of
these offsprings is the question of the maximal number of mutually
intersecting blocks of a $2-(v,k,1)$ design.

The Main Result above is contained in~\cite[Theorem\,6.4]{MR3336961};
the methods used there are different from our present methods. %
Our Main Result does not follow from~\cite[Chapter\,5]{MR3497070},
because the parameters $v=q^3+1$ and $k=q+1$ of unitals do not satisfy
the inequality $v > ((k-1)^2 +1)k$ in the assumption
of~\cite[5.3.5]{MR3497070}. %

\section{A class of linear spaces}
\label{sec:linearspaces}

A \emph{linear space} is an incidence geometry $(D,\cX,I)$ with point
set~$D$, line set~$\cX$ and incidence relation $I\subseteq D\times\cX$
such that any two members of~$D$ are incident with exactly one member
of~$\cX$, and each member of~$\cX$ is incident with at least two
members of~$D$. %

For $p\in D$ and $X\in\cX$, we denote the (line) pencil of~$p$ by
$\cX_p\coloneqq\set{L\in\cX}{(p,L)\in I}$, and the point row of~$X$ by
$D_X \coloneqq \set{a\in D}{(a,X)\in I}$. %
If the incidence relation is containment ``$\in$'' we identify lines
with their point rows, and we write $(D,\cX)$ instead of $(D,\cX,I)$.

\enlargethispage{5mm}%
\begin{assu}\label{assumptions}
  In this section, let $\DD \coloneqq (D,\cX,I)$ be a linear space with
  the following properties:
  \begin{itemize}
  \item $|D|=q^2$ for some integer $q>1$,
  \item For each point $p\in D$ we have $|\cX_p|\le q+1$,
  \item For each line $X\in\cX$ we have $|D_X|\le q+1$.
  \end{itemize}
\end{assu}

\begin{exas}\label{exasFromProjective}
  Let\/ $\PP=(P,\cL,J)$ be a projective plane of order~$q$. The
  following subgeometries (where $I\coloneqq (D\times\cX)\cap J$) are
  examples of linear spaces $(D,\cX,I)$ satisfying the assumptions
  in~\ref{assumptions}:
  \begin{enumerate}
  \item $D\coloneqq P\smallsetminus P_W$, $\cX\coloneqq\cL\smallsetminus\{W\}$,
    for some line $W\in\cL$.
  \item $D\coloneqq P\smallsetminus (P_W\cup\{v\}) \cup\{u\}$, %
    and $\cX=\cL$, where $W\in\cL$ and $u,v\in P$ such that $u\in P_W$
    but $v\notin P_W$.
  \item\label{ex:oval}%
    $D\coloneqq P\smallsetminus C$, $\cX\coloneqq\cL$, where~$C$ is a set
    of $q+1$ points in~$\PP$ such that no $q$ of them are collinear.
  \end{enumerate}
  As classical examples of sets~$C$ as
  in~\ref{exasFromProjective}.\ref{ex:oval}, we mention conics (or
  more generally, ovals) in projective planes.
\end{exas}

We will show that the examples in~\ref{exasFromProjective} cover all
possibilities. %
In other words, every $\DD$ as in~\ref{assumptions} is a subgeometry
of a projective plane of order~$q$.

\goodbreak
\begin{lemm}\label{lem:projectiveLine}
  If $L\in\cX$ satisfies $|D_L|=q+1$ then $D_L\cap D_X\ne\emptyset$
  holds for each $X\in\cX$.
\end{lemm}
\begin{proof}
  For $X\in\cX\smallsetminus\{L\}$ pick a point $p\in D_X\smallsetminus
  D_L$. Joining~$p$ to each point in~$D_L$ we obtain $q+1$ different
  members of~$\cX_p$. Now $|\cX_p|\le q+1$ yields that~$X$ is among
  these joining lines, and has a point in common with~$L$. 
\end{proof}

Every line with $q+1$ points will be called a \emph{projective line}
in~$\DD$.  This is justified by~\ref{lem:projectiveLine}. %

\begin{lemm}\label{thinPoint}
  There is at most one point\/ $u\in D$ with\/ $|\cX_u|\le q$. %
  If such a point\/~$u$ exists then $|\cX_u|=q$, there exists precisely
  one line $U\in\cX_u$ with\/ $|D_U|=q$, and all other lines through~$u$
  are projective lines.
\end{lemm}
\begin{proof}
  Assume that $s$ lines through~$u$ are not projective. Then the
  inequality $q^2 = |D| = 1 +\sum_{X\in\cX_u}|D_X\smallsetminus\{u\}| \le
  1+s(q-1)+(|\cX_u|-s)q \le 1+s(q-1)+(q-s)q = 1-s+q^2$ yields
  $s\le1$. As $s=0$ gives the contradiction $q^2=1+|\cX_u|q$, we have
  $s=1$. %
  Let~$S$ be the unique line in~$\cX_u$ with $|D_S|\le q$. %
  Now $q^2\le 1+(|D_S|-1)+(|\cX_u|-1)q = |D_S|+|\cX_u|q-q \le |\cX_u|q
  \le q^2$ yields $|\cX_u|=q$ and $|D_S|=q$.
\end{proof}

\begin{theo}\label{affinePlane}
  If\/~$\DD$ has no projective lines then~$\DD$ is an affine plane of
  order~$q$. In particular, we have $|\cX|=q^2+q$ in that case. 
\end{theo}
\begin{proof}
  From~\ref{thinPoint} we infer that $|\cX_p|=q+1$ for each $p\in D$
  if~$\DD$ has no projective lines. The inequality $q^2=|D| = 1
  +\sum_{X\in\cX_p}|D_X\smallsetminus\{p\}| \le 1+(q+1)(q-1) = q^2$ yields
  $|D_X|=q$ for each $X\in\cX$. Consider a point~$p$ and a line
  $X\in\cX\smallsetminus\cX_p$. Joining~$p$ with points on~$X$ we obtain
  all lines in~$\cX_p$ except for a single one. Thus the parallel
  axiom is satisfied, and~$\DD$ is an affine plane. 
\end{proof}

\begin{theo}\label{thinPencil}
  If there exists a point\/ $u\in D$ with $|\cX_u|\le q$ then~$\DD$ is
  obtained from a projective plane $\PP=(P,\cL)$ of order~$q$ by
  deleting some line~$W$ from~$\cL$ and the point set
  $(P_W\smallsetminus\{u\})\cup\{v\}$ from~$P$, where $u$ is some point
  on~$W$ and $v$ is some point in $P\smallsetminus P_W$. %
  In this case, we have $|\cX_u|=q$, and  $|\cX_p| = q+1$ holds for
  each $p\in D\smallsetminus\{u\}$.

  The line joining~$u$ and~$v$ has~$q$ points, every other line
  through~$u$ has $q+1$ points, and every other line through~$v$ has
  $q-1$ points. Each line in $\cX\smallsetminus(\cX_u\cup\cX_v)$ has~$q$
  points. 
\end{theo}
\begin{proof}
  Assume $|\cX_u|=q$. From~\ref{thinPoint} we know that there exists
  $S\in\cX_u$ with $|D_S|=q$, and that all other lines through~$u$ are
  projective. For $p\in D\smallsetminus D_S$, we obtain~$q$ lines
  in~$\cX_p$ by joining~$p$ to points on~$S$, and precisely one
  line~$A_p$ through~$p$ with $D_{A_p}\cap D_S=\emptyset$.  There
  is exactly one projective line in~$\cX_p$ (namely, the one
  joining~$p$ and~$u$) because any other such line would lead to
  $q+1$ joining lines in~$\cX_u$.

  We construct a new linear space $\Aa=(A,\cX)$, as follows. Pick a
  new symbol $v\notin D$ and form
  $A \coloneqq (D\smallsetminus\{u\})\cup\{v\}$. Retain the line
  set~$\cX$ and the old incidence relation on
  $({D\smallsetminus\{u\}})\times\cX$, and declare the new point~$v$
  incident with the members of %
  $\{S\}\cup\set{A_p}{p\in D\smallsetminus D_S}$ (and no others).

  Now every line in~$\Aa$ is incident with at most~$q$ points, and
  every pencil has at most $q+1$ lines. As $|A|=q^2$, we infer
  from~\ref{affinePlane} that~$\Aa$ is an affine plane of
  order~$q$. The original linear space~$\DD$ is obtained from the
  projective hull~$\PP$ of~$\Aa$ by deleting the line at infinity and
  all points at infinity apart from the one belonging to~$S$ (i.e.,
  the parallel class of~$S$ in~$\Aa$), and deleting the point~$v$.

  The numbers of points on each line and of lines through each point
  are now obvious.
\end{proof}

\begin{theo}\label{allFromProjective}
  Let $\DD=(D,\cX)$ be a linear space such that there is an integer
  $q\ge3$ with $|D|=q^2$, $|\cX_p|\le q+1$ for each $p\in D$, and
  $|D_X|\le q+1$ for each $X\in\cX$. %
  Then~$\DD$ is obtained by deleting $q+1$ points from a projective
  plane of order~$q$. Explicitly, we have the following cases:
  \begin{enumerate}
  \item\label{affine}%
    There are no projective lines, and~$\DD$ is an affine plane,
    with $|\cX|=q^2+q$.
  \item\label{thinPointItem}%
    There exists a point~$u$ with $|\cX_u|\le q$, and\/~$\DD$ is
    obtained as in~\ref{thinPencil}, with $|\cX|=q^2+q$. %
    In this case, we have $|\cX_u|=q$, and $|\cX_p| = q+1$ holds for
    each $p\in D\smallsetminus\{u\}$.
  \item\label{allLines}%
    We have $|\cX|=q^2+q+1$, and $|\cX_p|=q+1$ for each $p\in D$.  In
    this case, the linear space~$\DD$ is obtained from a projective
    plane of order~$q$ by deleting a set~$Y$ of $q+1$ points such that
    no~$q$ of them are collinear. %
    For each $y\in Y$, there exists $X_y\in\cX_y$ such that $|D_{X_y}|=q$. 
  \end{enumerate}
\end{theo}
\begin{proof}
  Assume that that~$\DD$ is not an affine plane. 
  From~\ref{affinePlane} we know that there exists a projective line
  $L\in\cX$, and obtain %
  $|\cX| = 1 + \sum_{p\in D_L}|\cX_p\smallsetminus\{L\}|$. %
  If there exists a point~$u$ with $|\cX_u|\le q$
  then~\ref{thinPencil} gives assertion~\ref{thinPointItem}. %

  There remains the case where $|\cX| = 1 + |D_L|\cdot q = 1+q^2+q$;
  we have $|\cX_p|=q+1$ for each $p\in D$ by~\ref{thinPencil}.  %
  So~$\DD$ is a $(q+1,1)$ design in the sense of~\cite{MR0345842},
  with precisely~$q^2$ points.  Then~$\DD$ is embeddable in a
  projective plane $\PP = (P,\cL)$ of order~$q$
  by~\cite[Theorem\,3.2]{MR0345842}; see
  also~\cite[Theorem\,3.1(i)]{MR0538056} (and~\cite{MR704236} for
  $q>3$). %

  Let~$Y \coloneqq P\smallsetminus D$ be the set of points of~$\PP$ that
  are not points of~$\DD$, and consider $y\in Y$. Joining~$y$ with the
  other~$q$ points in~$Y$, we obtain at most~$q$ of the $q+1$ lines
  in~$\cL_y$. Therefore, there exists a ``tangent'', i.e. a line
  $X_y\in\cX$ with $y\in P_{X_y}\smallsetminus D$. Now
  $|D_{X_y}\cap Y|=1$, and thus $|D_{X_y}|=q$.
\end{proof}

\begin{exam}\label{nearPencil}
  If $q=2$, the picture is slightly different from the one given
  in~\ref{allFromProjective}.

  In fact, if $\DD=(D,\cX,I)$ is a linear space such that $|D|=4$,
  $|\cX_p|\le3$ for each $p\in D$, and $|D_X|\le3$ for each
  $X\in\cX$ %
  then either $\DD$ is an affine plane of order~$2$, or it has one
  line of size~$3$ and three lines of size~$2$. %
  Note that this structure does not contain a quadrangle. 

  Except for the results about numbers of lines and sizes of line
  pencils (which fail because lines with less than two remaining
  points vanish from the picture), this structure is an example both
  for~\ref{allFromProjective}.\ref{thinPointItem} and
  for~\ref{allFromProjective}.\ref{allLines}. See
  Figure~\ref{fig:nearPencil}.
\end{exam}

\begin{rema}
  Doyen has determined all linear spaces with less than ten
  points. This is stated in~\cite{MR0249316}, but no details are
  given. Doyen's list is reproduced in~\cite[Appendix]{MR1253067}.
  According to that list, there exist\footnote{ \ %
    The pertinent examples are found on pages~211, 199, and~197
    in~\cite{MR1253067}. Note a tiny misprint in the example on
    page~197; the 13th line is missing in the picture.} %
  only two isomorphism types of linear spaces with~$9$ points and~$12$
  lines, and only one type with~$9$ points and~$13$ lines. %
  In the results of a computer-based search by Betten and
  Betten~\cite{MR1670277}, these three isomorphism types occur in
  Table\,IV as line case $(3^4,4^3)$ with~$13$ lines, and as cases
  $(3^7,4^2)$ and $(3^{12})$ with~$12$ lines, respectively.

  These examples are just the examples obtained (for $q=3$) in the
  three subcases of~\ref{allFromProjective}.  
\end{rema}

\section{Cliques in unitals}

Recall that a \emph{unital of order}~$q>1$ is a $2$-$(q^3+1,q+1,1)$
design, i.e., it has $q^3+1$ points, and $q+1$ points per block; the
number of blocks per point is then~$q^2$, and there are $q^2(q^2-q+1)$
blocks in total.

The \emph{confluence graph}~$\CG\UU$ of any incidence geometry
$\UU = (U,\cB)$ has~$\cB$ as its set of vertices, and two vertices are
adjacent if they are incident with a common element of~$U$. %
We consider a \emph{clique} in the confluence graph~$\CG\UU$ of a
unital $\UU=(U,\cB)$ of order~$q$; i.e., a set $\cC\subseteq \cB$ of
mutually intersecting blocks.

Our main result, as stated in the introduction, is part of the
following EKR result for unitals:

\begin{theo}\label{CliquesOfSizeQQ-EKR}
  Let\/~$\cC$ be a set of mutually intersecting blocks in a
  unital\/~$\UU$ of order~$q>2$, i.e., a clique in the confluence
  graph~$\CG\UU$.  Then $|\cC| \le q^2$, and\/ $|\cC| = q^2$ if and only
  if\/~$\cC$ is the pencil of all blocks through some point of the
  unital.
\end{theo}

\begin{proof}
  Recall Hoffman's ratio bound~\cite[Proposition\,1.1.7]{MR4350112}
  (see~\cite[Section\,2.4]{MR3497070}, \cite{MR4218548}) %
  for a strongly regular graph $\Gamma$ with (standard) parameters
  $(v,k,\lambda,\mu)$ and eigenvalues $k> r\geq s$: each clique $\cE$
  of $\Gamma$ has size at most $1+k/(-s)$, and if $\cE$ has exactly
  $1+k/(-s)$ vertices, then every vertex outside $\cE$ is adjacent to
  exactly $\mu/(-s)$ vertices of $\cE$. %

  The confluence graph~$\CG\UU$ of~$\UU = (U,\cB)$ is strongly
  regular, and its relevant parameters are $v={(q^2-q+1)q^2}$,
  $k={(q+1)^2(q-1)}$, $\mu=(q+1)^2$, $r=q^2-q-2$, and $s=-q-1$,
  see~\cite[8.5.4.A]{MR4350112}, where unitals of order~$q$ show up as
  Steiner $2$-designs $S(2,q+1,q^3+1)$.  %
  So Hoffman's ratio bound becomes
  $|\cC|\le1+{k}/{(-s)} = 1+(q+1)(q-1) = q^2$, %
  and if~$|\cC|$ attains that bound then~$\cC$ has the following
  property: %
  \[
    \text{every block in } \cB\smallsetminus\cC %
    \text{ meets precisely } \frac\mu{-s} = q+1 %
    \text{ blocks in } \cC . %
    \eqno{(*)}
  \]
  Aiming at a contradiction, we assume that~$\cC$ has size~$q^2$ and
  is not a pencil of blocks through a point in~$\UU$. %

  We consider the linear space $\DD = (\cC,\cX,\ni)$, where
  $\cX \coloneqq \set{x\in U}{2\le|\cB_x\cap\cC|}$ is the set of all
  points that are obtained as intersections of two elements of the
  clique, and the incidence relation ``$\ni$'' is containment
  ``$\in$'' reversed. %
  We abbreviate $\cC_x \coloneqq \cB_x\cap\cC$, and write %
  $\cX_B \coloneqq B\cap\cX$ for $B\in\cB$. %

  For each $x\in\cX$, there exists $B\in\cC\smallsetminus\cB_x$. As every
  element of~$\cC$ meets~$B$, we obtain
  $|\cC_x| \le |\cX_B| \le |B| = q+1$. Clearly, we have
  $|\cX_B| \le |B| = q+1$ for each $B\in\cC$. So~$\DD$ satisfies the
  standing assumptions of Section~\ref{sec:linearspaces},
  see~\ref{assumptions}. %

  By~\ref{allFromProjective}, we have one of the following cases:
  \begin{enumerate}
  \item[\ref{affine}]%
    If $(\cC,\cX,\ni)$ is an affine plane of order~$q$ then $|\cX_C|=q+1$
    and $|\cC_x|=q$ hold for each $C\in\cC$ and each $x\in\cX$. For
    any $L\in\cB_x\smallsetminus\cC$, property~$(*)$ yields that there is
    precisely one $B\in\cC\smallsetminus\cC_x$ that meets~$L$, and the
    intersection point does not belong to~$\cX$. This contradicts
    $|\cX_B|=q+1$.
  \item[\ref{thinPointItem}]%
    Let~$B$ be the unique element of~$\cC$ with $|\cX_B|=q$, and
    let~$y$ be the unique element of~$\cX_B$ with $|\cC_y|=q$. %
    Then $|\cX_C|=q+1$ holds for each $C\in\cC\smallsetminus\{B\}$. %
    Pick $L\in\cB_y\smallsetminus\cC$, then~$L$ meets the~$q$ members
    of~$\cC_y$ and one more element $C\in\cC\smallsetminus\cB_y$.  Now~$L$
    and~$C$ meet in a point $z\in\cX$ because $|\cX_C|=q+1$. Thus~$L$
    meets at least $|\cC_y|+|\cC_z| \ge q+2$ elements of~$\cC$. This
    contradicts~$(*)$.
    \goodbreak%
  \item[\ref{allLines}]%
    There remains the case where $|\cX_C| = q+1$ for each $C\in\cC$
    (recall that we assume $q>2$), but there exists $x\in\cX$ such
    that $|\cC_x|=q$. Pick any $B\in\cC\smallsetminus\cC_x$. %
    The~$q$ elements of~$\cC_x$ meet~$B$ in~$q$ of the $q+1$ points
    of~$\cX_B$. Joining the last remaining point on~$B$ with~$x$ we
    obtain a line $L\in\cB_x\smallsetminus\cC$. By~$(*)$, the line~$L$
    meets one further element~$C$ of~$\cC$ apart from those
    in~$\cC_x$. The intersection point belongs to~$\cX$ because
    $\cX_C$ contains all $q+1$ points. But this means that two
    elements of~$\cC$ pass through that intersection point, and~$L$
    meets at least $q+2$ elements of~$\cC$. This contradicts~$(*)$. %
    \qedhere
  \end{enumerate}
\end{proof}

\begin{coro}\label{noDualAffine}
  If\/ $\UU$ is a unital of order $q>2$ then~$\UU$ does not contain the
  dual of any linear space with~$q^2$ elements such that each line of
  that space has at most $q+1$ points. %

  In particular, the unital does not contain a dual affine plane of
  order~$q$, let alone a projective plane of order~$q$.
\end{coro}
\begin{proof}
  The lines of such a subgeometry would form a clique of order~$q^2$,
  and~\ref{CliquesOfSizeQQ-EKR} yields a contradiction. 
\end{proof}

\begin{theo}\label{reconstructUnital}
  Let\/ $\UU=(U,\cB)$ and\/~$\UU'=(U',\cB')$ be unitals of order~$q$
  and~$q'$, respectively. If\/~$q>2$ then every isomorphism from
  $\CG\UU$ to $\CG{\UU'}$ extends to an isomorphism from~$\UU$
  onto~$\UU'$.
\end{theo}
\begin{proof}
  Let $\beta\colon \CG\UU \to \CG{\UU'}$ be an isomorphism. Then
  $\beta$ is a bijection from~$\cB$ onto~$\cB'$, and $q=q'$ follows.
  For each point $u\in U$, the pencil~$\cB_u$ is a clique of
  size~$q^2$ in~$\CG\UU$, and~$\beta$ maps that pencil to a
  clique~$\cC_u$ of size~$q^2$
  in~$\CG{\UU'}$. From~\ref{CliquesOfSizeQQ-EKR} we know that there
  exists $u'\in U'$ such that $\cB'_{u'} = \cC_u$. Mapping $u$ to~$u'$
  is the point map of an isomorphism with block map~$\beta$, as
  required.
\end{proof}

\begin{coro}
  Any unital\/~$\UU=(U,\cB)$ of order~$q$ can be reconstructed, up to
  isomorphism, from its confluence graph~$\CG\UU$.
\end{coro}
\begin{proof}
  If $q>2$ then $\UU$ is isomorphic to $(\cM,\cB,\ni)$, where  $\cM$
  is the set of all cliques of size~$q^2$ in~$\CG\UU$. %
  If $q=2$ then $\CG\UU$ contains cliques of size~$q^2=4$ that are not
  pencils. However, the unitals of order~$2$ form a
  single isomorphism class, so~$\UU$ is determined by the sheer size
  of its confluence graph. 
\end{proof}

\begin{figure}
  \centering
  \begin{tikzpicture}[
      scale=1.5,%
          every node/.append style={circle, fill=white,
      draw=black, %
      inner sep=0pt, %
      minimum size=1em, %
    }]%

    \def\sd{sqrt(2)} %
    \def\lw{.7pt} %

    \begin{scope}[cm={1,0,0,1,(8,0)}]
    \node[coordinate] (P0) at (0,0) {0} ; %
    \node[coordinate] (P1) at ({\sd},0) {1} ; %
    \node[coordinate] (P2) at ({2*\sd},0) {2} ; %
    \node[coordinate] (P3) at ({3/2*\sd},{\sd*sqrt(3)/2}) {3} ; %
    \node[coordinate] (P4) at ({\sd},{\sd*sqrt(3)}) {4} ; %
    \node[coordinate] (P5) at ({\sd/2},{\sd*sqrt(3)/2}) {5} ; %
    \node[coordinate] (P6) at ({\sd},{\sd*sqrt(3)/3}) {6} ; %

    \draw[line width={2*\lw}, dashed] (P6) circle({\sd*1/3*sqrt(3)}) ; %

    \draw[line width={5*\lw}] (P0) -- (P1) -- (P2) ; %
    \draw[line width={3*\lw}, color=white] (P0) -- (P1) -- (P2) ; %
    \draw[line width={1*\lw}] (P0) -- (P1) -- (P2) ; %

    \draw (P2) -- (P3) -- (P4) ; %
    \draw (P2) -- (P6) -- (P5) ; %
    \draw (P4) -- (P5) -- (P0) ; %
    \draw (P0) -- (P6) -- (P3) ; %
    \draw (P1) -- (P6) -- (P4) ; %
    \draw[line width={3*\lw}]              (P0) -- (P6) ; %
    \draw[line width={1*\lw}, color=white] (P0) -- (P6) ; %
    \draw[line width={3*\lw}]              (P1) -- (P6) ; %
    \draw[line width={1*\lw}, color=white] (P1) -- (P6) ; %
    \draw[line width={3*\lw}]              (P2) -- (P6) ; %
    \draw[line width={1*\lw}, color=white] (P2) -- (P6) ; %

    \node[fill=black] (PP0) at (P0) {} ; %
    \node[fill=black] (PP1) at (P1) {\color{white}$u$} ; %
    \node[fill=black] (PP2) at (P2) {} ; %
    \node (PP3) at (P3) {} ; %
    \node (PP4) at (P4) {$v$} ; %
    \node (PP5) at (P5) {} ; %
    \node[fill=black] (PP6) at (P6) {} ; %
  \end{scope}
  
  \begin{scope}[cm={1,0,0,1,(14,0)}]
    every node/.append style={circle, fill=white,
      draw=black, %
      inner sep=0pt, %
    }]%
  
    \node[color=black] (1) at ( 0,0) {} ;%
    \node[color=black] (2) at ( 1,0) {} ;%
    \node[color=black] (3) at ( 2,0) {} ;%
    \node (4) at ( 0,1) {} ;%
    \node[color=black] (5) at ( 1,1) {} ;%
    \node (6) at ( 2,1) {} ;%
    \node (7) at ( 0,2) {} ;%
    \node (8) at ( 1,2) {} ;%
    \node (9) at ( 2,2) {} ;%
    \node (10) [coordinate] at ( 2.25,-0.5) {10} ; %
    \node (11) [coordinate] at (-0.5,-0.25) {11} ; %
    \node (12) [coordinate] at (-0.25, 2.5) {12} ; %
    \node (13) [coordinate] at ( 2.5, 2.25) {13} ; %
    \draw [line width={5*\lw}] (1) -- (2) -- (3) ;% 
    \draw [line width={3*\lw}, color=white] (1) -- (2) -- (3) ;% 
    \draw [line width={1*\lw}] (1) -- (2) -- (3) ;% 
    \draw [line width=1pt] (4) -- (5) -- (6) ;% 
    \draw [line width=1pt] (7) -- (8) -- (9) ;% 
    \draw [line width=1pt] (1) -- (5) -- (9) ;% 
    \draw [line width=1pt] (3) -- (5) -- (7) ;%
    %%----------------------------------------- 
    \draw [line width=1pt] (2) -- (4) ;%
    \draw [line width=1pt] (8) -- (6) ;%
    \draw [line width=1pt] (6) -- (2) ;%
    \draw [line width=1pt] (4) -- (8) ;%
    %%----------------------------------------- 
    \draw [line width=1pt] (1) -- (4) -- (7) ;%
    \draw [line width=1pt] (2) -- (5) -- (8) ;%
    \draw [line width=1pt] (3) -- (6) -- (9) ;%
    %%----------------------------------------- 
    \draw [line width=1pt, out=315, in= 45]  (6) to (10) [out=225, in=315] to (1) ;%
    \draw [line width=1pt, out=225, in=315]  (2) to (11) [out=135, in=225] to (7) ;%
    \draw [line width=1pt, out=135, in=225]  (4) to (12) [out= 45, in=135] to (9) ;%
    \draw [line width=1pt, out= 45, in=135]  (8) to (13) [out=315, in= 45] to (3) ;%
    %%----------------------------------------- 
    \draw [line width={3*\lw}]              (1) -- (5) ;% 
    \draw [line width={1*\lw}, color=white] (1) -- (5) ;% 
    \draw [line width={3*\lw}]              (2) -- (5) ;% 
    \draw [line width={1*\lw}, color=white] (2) -- (5) ;% 
    \draw [line width={3*\lw}]              (3) -- (5) ;% 
    \draw [line width={1*\lw}, color=white] (3) -- (5) ;% 
  \end{scope}

  \end{tikzpicture}
  \caption{The near pencil with~$4$ lines (joining the black points) inside a
    projective plane of order~$2$ (left, as
    in~\ref{allFromProjective}.\ref{thinPointItem}: the line~$W$ is
    dashed), and inside a unital of order~$2$ (right).}
  \label{fig:nearPencil}
\end{figure}
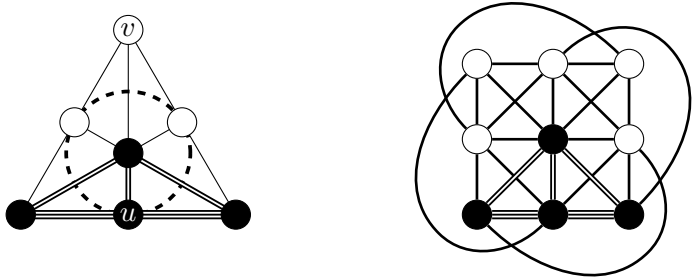

\begin{rema}
  The case $q=2$ forms a true exception in various respects. %
  First, we note that every unital of order~$2$ is isomorphic to the
  (unique) affine plane of order~$3$. %
  So it is trivially true that the isomorphism type of a given unital
  of order~$2$ is determined by its confluence graph and the
  information that we started with a unital (by the sheer number~$12$
  of vertices in the graph). However, not every automorphism of that
  graph extends to an automorphism of the unital: %

  The recognition of pencils in the confluence graph (as
  in~\ref{CliquesOfSizeQQ-EKR}) fails if $q=2$. 
  In fact, the confluence graph is then the complete multipartite
  graph $K_{3,3,3,3}$; its complementary graph is the disjoint union
  of~$4$ cliques of size~$3$ (namely, the parallel classes in the
  affine plane of order~$3$).

  Although the first assertion of~\ref{noDualAffine} fails for $q=2$
  (see~\ref{nearPencil} and Figure~\ref{fig:nearPencil}), the unital of
  order~$2$ does not contain any dual affine plane of order~$2$. In
  fact, the latter structure is a set of~$6$ points such that no three of
  them are collinear. In an affine plane of order~$3$, however, any
  set of at least~$5$ points contains three collinear points.
\end{rema}

\section{Partial linear spaces without O'Nan configurations}%%

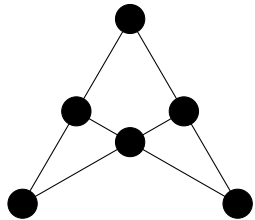
\begin{figure}[h!]
  \centering
  \begin{tikzpicture}[
      scale=1.2,%
          every node/.append style={circle, fill=black,
      draw=black, %
      inner sep=0pt, %
      minimum size=1em, %
    }]%

    \def\sd{sqrt(2)} %
    \def\lw{.7pt} %

    \node[coordinate] (P0) at (0,0) {0} ; %
    \node[coordinate] (P2) at ({2*\sd},0) {2} ; %
    \node[coordinate] (P3) at ({3/2*\sd},{\sd*sqrt(3)/2}) {3} ; %
    \node[coordinate] (P4) at ({\sd},{\sd*sqrt(3)}) {4} ; %
    \node[coordinate] (P5) at ({\sd/2},{\sd*sqrt(3)/2}) {5} ; %
    \node[coordinate] (P6) at ({\sd},{\sd*sqrt(3)/3}) {6} ; %

    \draw (P2) -- (P3) -- (P4) ; %
    \draw (P2) -- (P6) -- (P5) ; %
    \draw (P4) -- (P5) -- (P0) ; %
    \draw (P0) -- (P6) -- (P3) ; %

    \node (PP0) at (P0) {} ; %
    \node (PP2) at (P2) {} ; %
    \node (PP3) at (P3) {} ; %
    \node (PP4) at (P4) {} ; %
    \node (PP5) at (P5) {} ; %
    \node (PP6) at (P6) {} ; %
  \end{tikzpicture}
  \caption{The O'Nan configuration.}
  \label{fig:O'Nan}
\end{figure}

Recall that the O'Nan configuration has~$6$ points and~$4$ mutually intersecting
lines, such that each line is incident with exactly~$3$ points, and
each point is incident with exactly~$2$ lines. See
Figure~\ref{fig:O'Nan}. 

A \emph{partial linear space} is an incidence geometry $(D,\cX,I)$
where any two members of~$D$ are incident with \emph{at most} one
member of~$\cX$, and each member of~$\cX$ is incident with at least
two members of~$D$.  For $(p,L)\in (D\times\cX)\smallsetminus I$, the
\emph{near pencil determined by} $(p,L)$ consists of the line~$L$, and
all lines joining points on~$L$ with~$p$.  % 

\goodbreak%
\begin{prop}\label{noONan}
  Let\/~$(P,\cL)$ be a partial linear space, and let $\cC\subseteq\cL$
  be a set of mutually intersecting lines.
  \begin{enumerate}
  \item If\/~$(P,\cL)$ contains no triangle then~$\cC$ is contained in
    a pencil; thus the maximal cliques in the confluence graph
    of\/~$(P,\cL)$ are the pencils. This applies to all generalized
    $m$-gons\/~$(P,\cL)$ with $m>3$, and to all near $m$-gons with
    $m>3$.
  \item If\/\/~$(P,\cL)$ contains a triangle, but no O'Nan
    configuration then~$\cC$ is contained in a pencil, or in a near
    pencil comprising at least~$3$ lines.
  \end{enumerate}
\end{prop}
\begin{proof}
  Without loss, we may assume that~$\cC$ has more than one element. %
  Pick distinct $K,L\in\cC$ and let~$p$ be the unique point in $K\cap L$. 
  If~$\cC$ does not contain a triangle then each member of~$\cC$
  passes through~$p$. %
  If~$\cC$ contains a triangle then the absence of O'Nan
  configurations implies that at most one of the vertices of
  that triangle is on more than two members of~$\cC$. So~$\cC$ is
  contained in a near pencil.
\end{proof}

\begin{coro}\label{maxCliquesHermitian}
  Let\/~$\UU$ be a unital of order~$q$ without O'Nan configurations
  (e.g. a Hermitian unital over a finite field). %
  Then the maximal cliques in the confluence graph~$\CG\UU$ of\/~$\UU$
  are the pencils and the near pencils. %
\end{coro}

\begin{proof}
  In a Hermitian unital (over any commutative field) there are no
  O'Nan configurations, see~\cite{MR0295934}, cp.~\cite{MR2795696}. %
  Every near pencil contains $q+2$ blocks, but at most $q+1$ of those
  blocks belong to any single pencil.  Since the size $q^2$ of a
  pencil is greater than $q+1$, no pencil is contained in a near
  pencil. Thus~\ref{noONan} yields the assertion.
\end{proof}

Our description above of the maximal cliques of Hermitian unitals
becomes false for unitals with O'Nan configurations, because the
(maximal) cliques containing such a configuration are neither
contained in a pencil nor in a near pencil.

It has been conjectured (cf.~\cite[p.\,102]{MR587626},
\cite[p.\,87]{MR2440325}) that finite Hermitian unitals are
characterized by the absence of O'Nan configurations. To our
knowledge, this conjecture has not been proved yet. %

%%%%%%%%%%%%%%%%%%%%%%%%%%%%%%%%%%%%%%%%%%%%%%%%%%%%%%%%%%%%%%%%%%

%%%%%%%%%%%%%%%%%%%%%%%%%%%%%%%%%%%%%%%%%%%%%%%%%%%%%%%%%%%%%%%%%%%%%%%
\vfill

\begin{tabular}{*{3}{p{.3\textwidth}}}
{Theo Grundh\"ofer}&%
{Markus Stroppel}&%
{Hendrik Van Maldeghem}\\
Institut f\"ur Mathematik&%
LExMath &%
          Vakgroep Wiskunde: \\
                   & %
                     Fakult\"at 08 &%
  Algebra en Meetkunde\\
Universit\"at W\"urzburg&%
Universit\"at Stuttgart   &%
Universiteit Gent\\
  Am Hubland&%
              &%
Krijgslaan 281, S25\\
D-97074 W\"urzburg&%
D-70550 Stuttgart&%
B--9000 Gent\\
Germany&%
Germany&%
Belgium
\end{tabular}

%%%%%%%%%%%%%%%%%%%%%%%%%%%%%%%%%%%%%%%%%%%%%%%%%%%%%%%%%%%%%%%%%%%%%%% 

\begin{thebibliography}{10}
\providecommand{\href}[2]{#2}
\providecommand{\eprint}[1]{\href{http://arxiv.org/abs/#1}{#1}}
\providecommand{\doi}[1]{\url{https://doi.org/#1}}
% we adapt \MR , making it work well with MathSciNet bibtex version:
\providecommand{\MR}[1]{\relax\ifhmode\unskip\space\fi \MRnumberextract#1 \,}
\def\MRnumberextract#1 #2\,{\MRhref{#1}{#2}}%
\providecommand{\MRhref}[2]{%
  \href{https://mathscinet.ams.org/mathscinet-getitem?mr=#1}{MR\,#1 #2}}
\providecommand{\ZBL}[1]{\relax\ifhmode\unskip\space\fi \ZBLhref{#1}}
\providecommand{\ZBLhref}[1]{%
  \href{http://zbmath.org/?q=an:#1}{Zbl #1}}
\providecommand{\JfM}[1]{\relax\ifhmode\unskip\space\fi \JfMhref{#1}}
\providecommand{\JfMhref}[1]{%
  % \href{https://www.emis.de/cgi-bin/jfmen/MATH/JFM/quick.html?type=html&an=#1}{JfM #1}
  \href{http://zbmath.org/?q=an:#1}{JfM #1}}

\bibitem{MR2440325}
S.~G. Barwick and G.~Ebert, \emph{Unitals in projective planes}, Springer
  Monographs in Mathematics, Springer, New York, 2008,
  \doi{10.1007/978-0-387-76366-8}. \MR{2440325.} \ZBL{1156.51006}.

\bibitem{MR1253067}
L.~M. Batten and A.~Beutelspacher, \emph{The theory of finite linear
  spaces}, Cambridge University Press, Cambridge, 1993,
  \doi{10.1017/CBO9780511666919}. \MR{1253067.} \ZBL{0806.51001}.

\bibitem{MR1670277}
A.~Betten and D.~Betten, \emph{Linear spaces with at most 12 points}, J.
  Combin. Des. \textbf{7} (1999), no.~2, 119--145, \\
  \doi{10.1002/(SICI)1520-6610(1999)7:2<119::AID-JCD5>3.0.CO;2-W}.
  \MR{1670277.} \ZBL{0935.51009}.

\bibitem{MR4350112}
A.~E. Brouwer and H.~Van~Maldeghem, \emph{Strongly regular graphs},
  Encyclopedia of Mathematics and its Applications  182, Cambridge University
  Press, Cambridge, 2022, \doi{10.1017/9781009057226}. \MR{4350112.}
  \ZBL{7437385}.

\bibitem{MR3336961}
M.~De~Boeck, \emph{The largest {E}rd{\H{o}}s-{K}o-{R}ado sets in 2-{$(v,k,1)$}
  designs}, Des. Codes Cryptogr. \textbf{75} (2015), no{}~3, 465--481,
  \doi{10.1007/s10623-014-9929-5}. \MR{3336961.} \ZBL{1312.05022}.

\bibitem{MR704236}
S.~Dow, \emph{An improved bound for extending partial projective planes},
  Discrete Math. \textbf{45} (1983), no. 2-3, 199--207,
  \doi{10.1016/0012-365X(83)90036-5}. \MR{704236.} \ZBL{0511.05021}.

\bibitem{MR0249316}
J.~Doyen, \emph{Sur le nombre d'espaces lin\'eaires non isomorphes de {$n$}\
  points}, Bull. Soc. Math. Belg. \textbf{19} (1967), 421--437. \MR{0249316.}
  \ZBL{0157.03401}.

\bibitem{MR3497070}
C.~Godsil and K.~Meagher, \emph{Erd\H{o}s-{K}o-{R}ado theorems: algebraic
  approaches}, Cambridge Studies in Advanced Mathematics  149, Cambridge
  University Press, Cambridge, 2016, \doi{10.1017/CBO9781316414958}.
  \MR{3497070.} \ZBL{1343.05002}.

\bibitem{MR2795696}
T.~Grundh{\"o}fer, B.~Krinn, and M.~J. Stroppel, \emph{Non-existence of
  isomorphisms between certain unitals}, Des. Codes Cryptogr. \textbf{60}
  (2011), no.~2, 197--201, \doi{10.1007/s10623-010-9428-2}. \MR{2795696.}
  \ZBL{05909195}.

\bibitem{GrundhoeferStroppelVanMaldeghem-confluence-arxiv-v1}
T.~Grundh{\"o}fer, M.~J. Stroppel, and H.~{Van Maldeghem},
  \emph{Confluence graphs of unitals}, 2023. arXiv~\eprint{2311.11693}.

\bibitem{MR4218548}
W.~H. Haemers, \emph{Hoffman's ratio bound}, Linear Algebra Appl. \textbf{617}
  (2021), 215--219, \doi{10.1016/j.laa.2021.02.010}. \MR{4218548.}
  \ZBL{1459.05172}.

\bibitem{MR0538056}
D.~McCarthy, N.~M. Singhi, and S.~A. Vanstone, \emph{A graph-theoretical
  approach to embedding {$(r,\,1)$}-designs}, in: \emph{Graph theory and
  related topics ({P}roc. {C}onf., {U}niv. {W}aterloo, {W}aterloo, {O}nt.,
  1977)}, pp. 289--304, Academic Press, New York-London, 1979.
  \MR{538056.} \ZBL{0459.05016}.

\bibitem{MR0295934}
M.~E. O'Nan, \emph{Automorphisms of unitary block designs}, J. Algebra
  \textbf{20} (1972), 495--511, \doi{10.1016/0021-8693(72)90070-1}.
  \MR{0295934.} \ZBL{0241.05013}.

\bibitem{MR587626}
F.~Piper, \emph{Unitary block designs}, in: \emph{Graph theory and
  combinatorics ({P}roc. {C}onf., {O}pen {U}niv., {M}ilton {K}eynes, 1978)},
  Res. Notes in Math. 34, no. 98--105, Pitman, Boston, Mass., 1979.
  \MR{587626.} \ZBL{0455.05024}.

\bibitem{MR3135960}
B.~Stroppel and M.~J. Stroppel, \emph{Desargues, doily, dualities and
  exceptional isomorphisms}, Australas. J. Combin. \textbf{57} (2013),
  251--270,
  \url{https://ajc.maths.uq.edu.au/pdf/57/ajc_v57_p251.pdf}.
  \MR{3135960.} \ZBL{06238433}.

\bibitem{MR0345842}
S.~A. Vanstone, \emph{The extendibility of {$(r,\,1)$}-designs}, in:
  \emph{Proceedings of the {T}hird {M}anitoba {C}onference on {N}umerical
  {M}athematics ({W}innipeg, {M}an., 1973)}, Congress. Numer.~IX, pp.
  409--418, Utilitas Math., Winnipeg, MB, 1974. \MR{345842.} \ZBL{0317.05015}.

\end{thebibliography}
\end{document}